# Akademie Intakt 2021
## English Edition

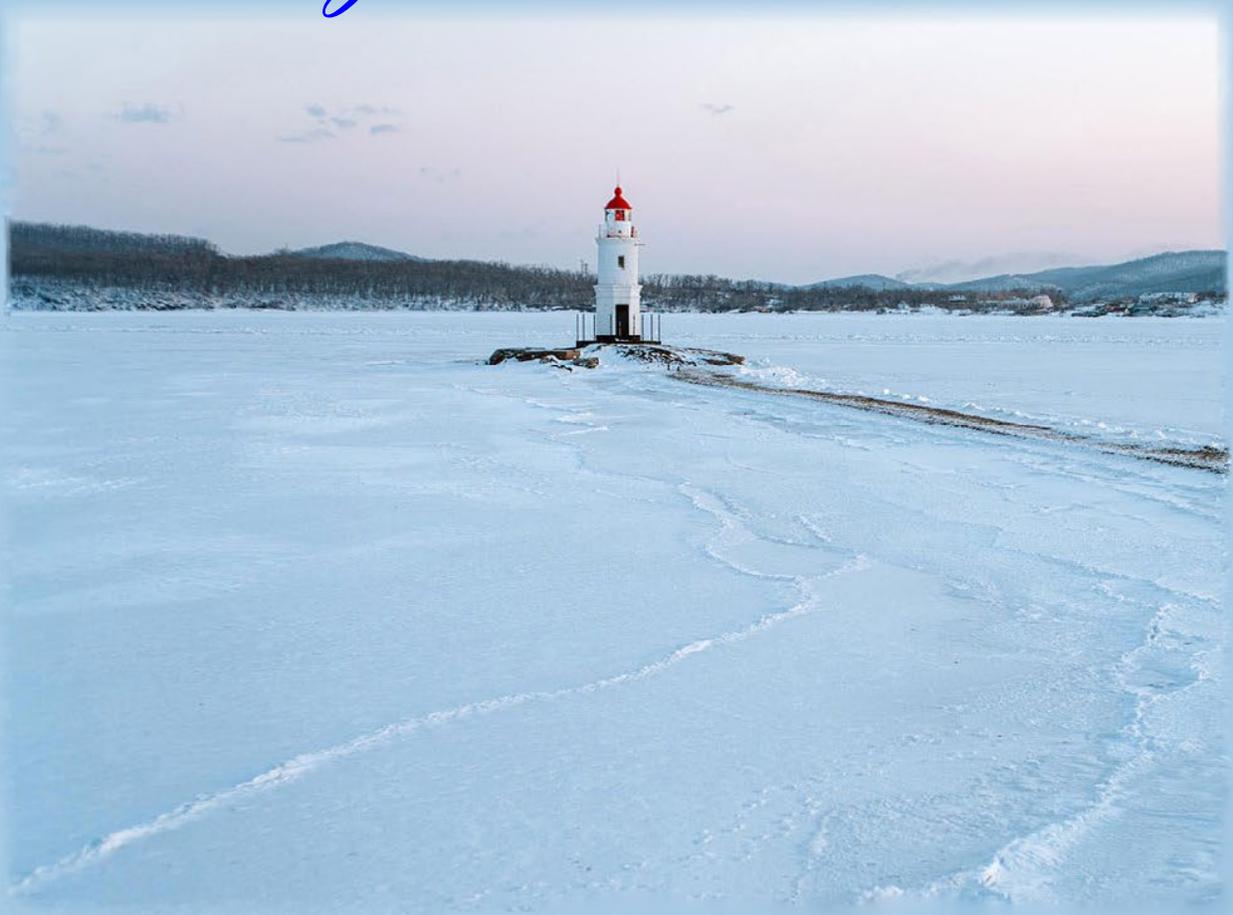

# Celebrating Women in Mathematics

Since 2019, 12 May is the day for celebrating women in mathematics. Why 12 May and how did this initiative appear?

12 May is the birthdate of Maryam Mirzakhani, the first female mathematician who was awarded with the prestigious Fields Medal. The Fields Medal is granted every four years at the International Congress of the International Mathematical Union (IMU) to scientists under 40 years for outstanding contributions in mathematics. In the absence of a Nobel Prize in mathematics, the Fields Medal is considered to be the most prestigious international award in the field of mathematics.

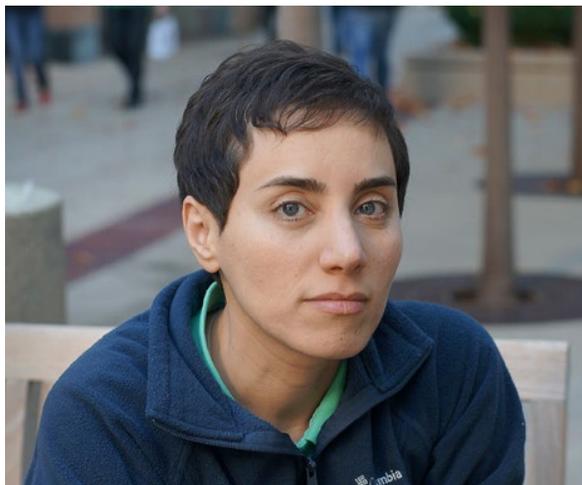
Maryam Mirzakhani © Stanford News Service

Maryam Mirzakhani (1977–2017) (Fig. 1) was born in Tehran, Iran, and grew up there. While in high school, she won twice gold medals for Iran in the International Mathematical Olympiad (IMO): at IMO 1994 in Hong Kong, she was the first Iranian female who won a gold medal with 41 points out of the maximum of 42; at IMO 1995 in Toronto, she won a gold medal gaining the maximum score of 42 points and became the first Iranian with two gold medals from IMOs. Mirzakhani continued her education at the Sharif University of Technology in Tehran, where she completed a bachelor's degree in mathematics in 1999. Then she moved to the United States, completed a PhD in 2004 at Harvard University, and got hired as an assistant professor at Princeton University. She was also appointed as a Clay Research Fellow. In 2009, Mirzakhani became a professor at Stanford University and in 2013 she was awarded by the Simons Foundation to be a Simons Investigator.

Mirzakhani was honored with the Fields Medal at the age of 37 "for her outstanding contributions to the dynamics and geometry of Riemann surfaces and their moduli spaces" (1). This happened on 13 August 2014 at the Opening Ceremony of the International Congress of Mathematicians (ICM) in Seoul, South Korea (Fig. 2). Three years later Mirzakhani died from cancer and the world lost one of the greatest mathematicians.

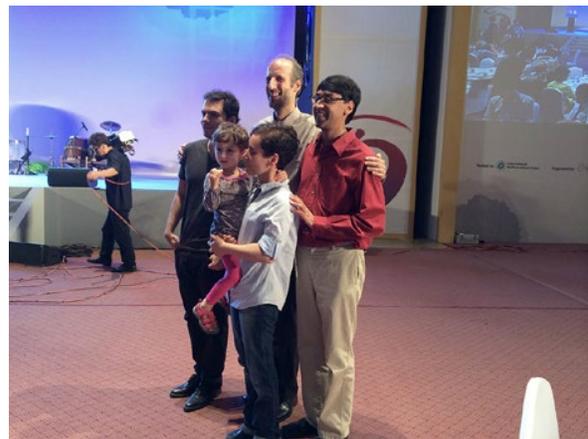
Maryam Mirzakhani among the recipients of the Fields Medal at the 2014 International Congress of Mathematicians in Seoul © Wikipedia/Monsoon0

During the World Meeting for Women in Mathematics (WM) (2) on 31 July 2018 (a satellite event of the ICM 2018), some activities were set in memory of Mirzakhani. The IMU Committee for Women in Mathematics created the exhibition "Remember Maryam Mirzakhani" with 18 posters portraying Mirzakhani, volumes with her mathematical papers, as well as a book with articles about her. The exhibition was opened at the WM (2) and remained open during the ICM 2018. As another tribute to the memory of Mirzakhani, the Women's Committee of the Iranian Mathematical Society presented a proposal for celebrating women in mathematics worldwide on Mirzakhani's birthdate (12 May) and this was approved by a large majority of the attendees.

This is how the May12 initiative arose. It became an international annual initiative (2), coordinated and supported by several organizations for woman in mathematics. Already in the first year, in 2019, there were more than 100 events worldwide. Despite the Covid-19 pandemic, the initiative continued in 2020 and 2021, though mostly with online activities. For the second edition in 2020, 152 events were published on the May12 website and more than 20,000 people from 131 countries expressed interest to the online screenings of "Secrets of the Surface: The Mathematical Vision of Maryam Mirzakhani" (a documentary film by George Csicsery on the life and work of Mirzakhani). For the year 2021, there are also numerous events listed on the May12-website, as well as screenings of the movie "Picture a Scientist" (by Sharon Shattuck and Ian Cheney) at which three female scientists (biologist, chemist, and geologist) share experiences from their journey in science.

The Austrian Academy of Sciences (OeAW) was also involved in the third edition of the May12 initiative. Dedicated to the International Day of Women in Mathematics and Maryam Mirzakhani, in cooperation with colleagues from the University of Vienna and the University of Novi Sad, we organized the Generalized Functions Online Workshop (3) which took place on the 12th of May 2021. The aim of the workshop was to increase the visibility of the female members of the Generalized Functions (GF) community and to facilitate scientific collaboration and networking within this group. It was also aimed to encourage early career female mathematicians, to give them the chance to present their first results, and to encourage them to pursue a scientific career in mathematics. A main focus of the conference was the panel discussion `"How to be a Good Ally to Women in Mathematics"'. However, participation was not limited to female members of the GF community. The workshop attracted the attention of scientists from many countries. There were participants from Austria, Belgium, Bosnia and Herzegovina, China, Kazakhstan, India, Italy, Poland, Romania, Serbia, the United Kingdom, the United States, and many other countries. We tried to make the online workshop as similar as possible to an in-person conference. Between the talks sessions, there were possibilities for interaction of the participants in smaller groups. We also organized virtual excursions and some math-fun activities for the participants. Although an online conference could not fully substitute in-person meetings and the great impact they have on scientific collaboration, it made possible presentations and interactions between scientists from different places of the world during times when traveling was not possible due to the Covid-19 pandemic. Even more, the online way actually made possible the participation of scientists who usually do not have financial support for conference travels (which very much concerns early career scientists and scientists from low and middle-income countries) or cannot travel due to family or other duties. Thus, we realized that even when traveling becomes possible again and we are back to the standard way, we should try to organize future events in a hybrid way in order to give the chance for participation to more people.

After paying tribute to Maryam Mirzakhani, we will continue with brief information on several other world-renowned female mathematicians, through a short walk from ancient times to the present day.

**Hypatia** (c. 350/370–415 AD) was a mathematician, philosopher, and astronomer from Alexandria, Egypt (province of the Roman Empire at that time). She was tutored in mathematics by her father, the mathematician Theon of Alexandria, who took care to educate very well his daughter in contrast to the traditional way women were raised at that time without much education. Hypatia's possibly most used portrait is shown in Fig. 3, which however is regarded as being a fictionalized sketch.

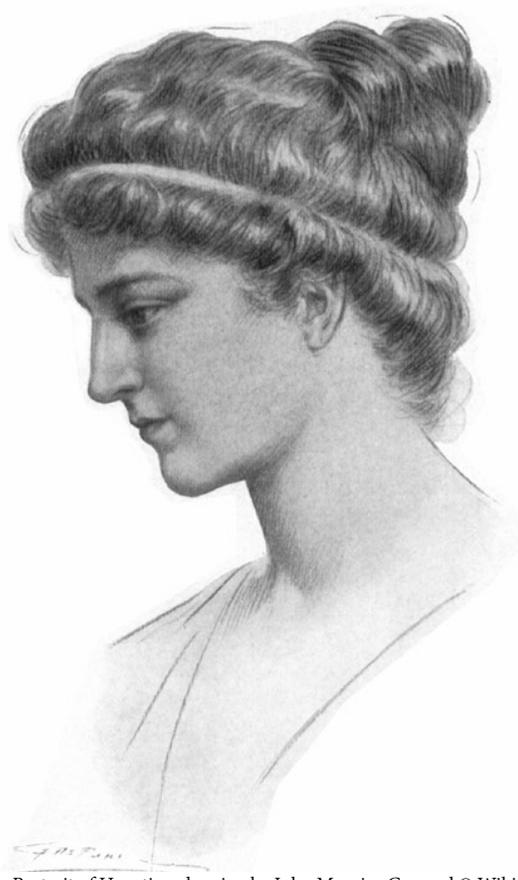

Portrait of Hypatia—drawing by Jules Maurice Gaspard © Wikipedia

Hypatia is considered to be the earliest known female mathematician, who contributed substantially to the development of mathematics and whose life and work are relatively well documented in ancient writing. According to (4), Hypatia "is credited with commentaries on Apollonius of Perga's Conics (geometry) and Diophantus of Alexandria's Arithmetic (number theory), as well as an astronomical table (possibly a revised version of Book III of her father's commentary on the Almagest)". Hypatia was also a respected popular teacher and lecturer on philosophical topics, attracting many students and large audience (Fig. 4, next page). It is believed she was the first woman who has taught mathematics.

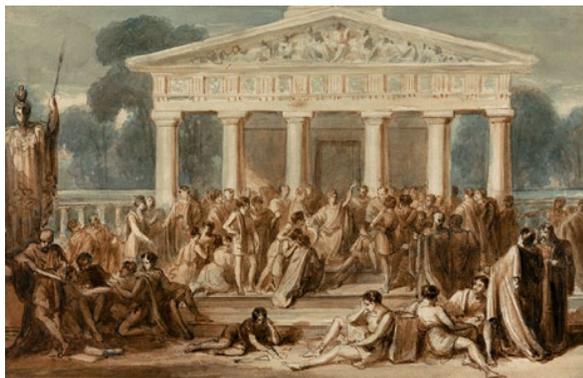
Hypatia teaching in Alexandria—painting by Robert Trewick Bone
© Britannica.com

**Sofia Kovalevskaya** (1850–1891) (Fig. 5) was a mathematician from Russia. Her parents took care to provide her with a good early education and Kovalevskaya showed to be talented for mathematics already in her childhood. However, at that time, women were not admitted at universities in some countries, including Russia. In order to be able to continue her education, Kovalevskaya married in 1868 and went to Germany to study mathematics there. In Germany there were also restrictions for women at some universities.

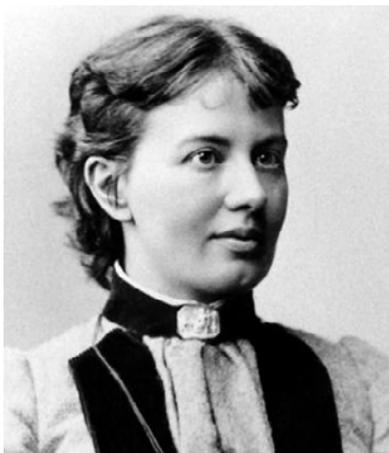
Sofia Kovalevskaja © Wikipedia

Kovalevskaya was not allowed to attend the University of Berlin, but she did not let this destroy her passion for math and her willingness for in-depth study. She proceeded studying mathematics privately under the guidance of the German mathematician Karl Weierstrass (one of the founders of the modern theory of functions and one of the most renowned mathematicians of his time). Under the supervision of Weierstrass and his encouragement, Kovalevskaya made a doctorate. Instead of the requirements to complete one dissertation, she actually completed three works during the 2 years' residence in Berlin, two in pure mathematics and one on a topic in astronomy. Based on these works, the University of Göttingen exempted her from the required examination and public defense of the dissertation and awarded to her the degree Doctor of Philosophy in 1874. In 1883, Kovalevskaya was invited to Stockholm to lecture mathematics. In 1889 she was appointed professorship in mathematics at the Stockholm University College (now Stockholm University).

Kovalevskaya made significant contributions to the theory of differential equations. She is regarded as being the greatest female mathematician prior to the twentieth century. She was the first woman who obtained a doctorate in mathematics, the first woman in the world who was appointed full professor in mathematics, the first female professor in Sweden, and one of the first women serving as editor of a scientific journal ("Acta Mathematica").

**Emmy Noether** (1882–1935) (Fig. 6) was a German mathematician and physicist. Her father was a professor of mathematics. Noether got certified to teach French and English in 1900, but she decided to continue with studies in mathematics. Though German universities were not fully opened for women yet, the situation was gradually improving at that time and women could attend lectures under the specific permission of the teaching professor.

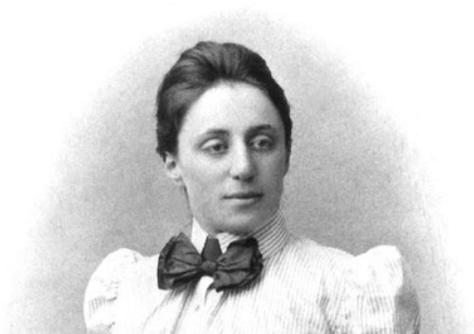
Emmy Noether © Wikipedia

In this way, Noether attended some lectures at the Universities in Erlangen and Göttingen. With the change of the rules, Noether was officially enrolled to study mathematics at the University of Erlangen and in 1907 she received a Ph.D. degree with a dissertation on algebraic invariants. Even with her doctorate, the University of Erlangen would not hire Noether, because of its policy against female professors. So, for the next several years, Noether worked unpaid in Erlangen, helping her father with teaching and working on her own research.

In 1915, she was invited by the renowned mathematicians David Hilbert and Felix Klein to join the University of Göttingen and to work with them on the mathematics behind Albert Einstein's theory of general relativity. There were some issues regarding the failure of local energy conservation in the general theory, and Hilbert and Klein thought her expertise on invariant theory would be helpful. Noether joined them and not only resolved the issue with Einstein's theory, but also obtained more general theorems. She discovered a link between symmetries and conservation laws (5), which is a key result in theoretical physics used by physicists ever since and known as Noether's theorem.

Despite the improvement of possibilities for women to study at universities, female mathematicians were still not very well accepted for habilitation and academic positions. Because of that Noether spent some years in Göttingen without a fixed university position.

Her habilitation was approved in 1919 and finally she was hired in 1923 at the University of Göttingen and she was able to officially guide students. In 1933, when the Nazi party started to rule Germany, Noether and many other professors were dismissed from the university and Noether accepted a visiting research position at Bryn Mawr College in Pennsylvania, USA. In 1935, at the age of 53, she died several days after an operation of a tumor.

Noether is mostly known for her important contributions to Abstract Algebra and her fundamental result in mathematical physics, the Noether's theorem. In 1932, she was awarded the prestigious Ackermann-Teubner Memorial Prize. For her innovative work on algebra, Noether is considered to be "the most creative abstract algebraist of modern times". She is also the first woman who gave a plenary lecture at the ICM (in 1932 in Zurich, Switzerland). Albert Einstein wrote in the New York Times in appreciation of Noether (6): "In the judgment of the most competent living mathematicians, Fräulein Noether was the most significant creative mathematical genius thus far produced since higher education of women began."

**Karen Keskulla Uhlenbeck** (born in 1942) is an American mathematician. She is a Professor and Sid W. Richardson Regents Chairholder at the Department of Mathematics, University of Texas, Austin, and a current Distinguished Visiting Professor at the Institute for Advanced Study in Princeton.

Uhlenbeck is one of the founders of the field of geometric analysis. In 2019, she became the first woman awarded with the prestigious Abel Prize "for her pioneering achievements in geometric partial differential equations, gauge theory and integrable systems, and for the fundamental impact of her work on analysis, geometry and mathematical physics" (7). She was awarded with the Emmy Noether Lecture in 1988. Uhlenbeck is the second female plenary speaker at the ICM (in 1990 in Kyoto, Japan), after Emmy Noether, who was the first one with such an honour in 1932. The ICM is the largest mathematical event, dating back to 1897 when the first edition took place in Zürich. Having female plenary speakers for the first time in 1932 and the second time in 1990 indicates how difficult the acceptance of women in a male-dominated field was.

**Ingrid Daubechies** (born in 1954) is a Belgian mathematician and physicist. Currently, she is a James B. Duke Distinguished Professor of Mathematics and Electrical and Computer Engineering at Duke University (North Carolina, USA). She is a member of the National Academy of Sciences, the National Academy of Engineering, and the American Academy of Arts and Sciences.

Daubechies is mostly known with her works on wavelets and their applications, and she is one of the world's most cited mathematicians. Daubechies is the first woman who became full professor in mathematics at Princeton University (in 1994), the first female elected for President of the International Mathematical Union (in 2011) and the first woman who received the National Academy of Sciences Award in Mathematics (in 2000). She was awarded the Emmy Noether Lecture in 2006, and numerous other prizes and honors.

The aforementioned women are just a small sample of all great female mathematicians. Their contribution to the field of mathematics is very significant; it will long be remembered and will have a great impact on future work in mathematics. With their dedication to mathematics and their great achievements, these successful female scientists can inspire and encourage other women to pursue mathematics with passion and to reach for the top.

**Diana T. Stoeva, Acoustics Research Institute of OeAW**